%% file: BEM.tex
\RequirePackage{fix-cm}
\documentclass[smallextended]{svjour3}       
\smartqed  
\usepackage{graphicx}
\usepackage{xcolor}
\usepackage{amsmath,amssymb}
\usepackage[english]{babel}
\usepackage[utf8]{inputenc}

\usepackage{diagbox}
\usepackage{multirow}
\usepackage{lscape}
\usepackage{float}

\usepackage{vmargin}
\setmarginsrb{70pt}{70pt}{70pt}{70pt}{20pt}{20pt}{20pt}{20pt}
\include{customdefs}
\begin{document}

\title{Least squares moment identification of binary regression mixtures models
}



\author{Benjamin Auder  \and Elisabeth Gassiat \and Mor Absa Loum}

\institute{\textbf{Benjamin Auder}\at
  \email{benjamin.auder@u-psud.fr}   \\
 Université Paris-Saclay, CNRS, Laboratoire de Mathématiques d'Orsay, 91405, Orsay , France
\and
  \textbf{Elisabeth Gassiat} \at
  \email{elisabeth.gassiat@u-psud.fr}\\
  Université Paris-Saclay, CNRS, Laboratoire de Mathématiques d'Orsay, 91405, Orsay , France
\and
\textbf{Mor Absa Loum} (corresponding author)  \at
  \email{morabsa.loum5@gmail.com} \\
  Université Paris-Saclay, CNRS, Laboratoire de Mathématiques d'Orsay, 91405, Orsay , France
}

\date{Received: date / Accepted: date}

\maketitle

\begin{abstract}
We consider finite mixtures of generalized linear models with binary output.
We prove that cross moments (between the output and the regression variables) until order 3 are sufficient to identify all parameters of the model.
We propose a least-squares estimation method based on those moments and we prove the consistency and the Gaussian asymptotic behavior of the estimator. 
We provide simulation results and comparisons with likelihood methods. Numerical experiments were conducted using the R-package morpheus that we developed
for our least-squares moment method and with the R-package flexmix for likelihood methods.
We then give some possible extensions to finite mixtures of regressions with binary output including both continuous and categorical covariates, and possibly longitudinal data.

\end{abstract}

\newpage
\input{intro}

\input{identification}

\input{simulations}

\input{extensions}

\input{proofs}

%
%

\bibliographystyle{spmpsci}      
\bibliography{biblio}   

%
%





\end{document}

%% file: customdefs.tex
\newcommand{\E}{\mathrm{E}}

%% file: intro.tex
\section{Introduction}
\label{intro}
Logistic models, or  more generally multinomial regression models that fit covariates to discrete responses through a link function,
are very popular for use in various application fields. When the data under study come from several latent groups that have different characteristics,
using mixture models is also a very popular way to handle heterogeneity.
Thus, many algorithms were developed to deal with various mixture models, see for instance the books \cite{MR2265601} and \cite{MR3889980}.
Most of them use likelihood or Bayesian methods that are likelihood dependent.
Indeed, the well known  expectation-maximization (EM) methodology or its randomized versions makes it often easy to build algorithms.
However, one problem of such methods is that they can converge to local  spurious maxima so that it is
necessary to explore multiple initial points as explained in \cite{WU1983}. This is one reason for
renewed interest in moment based methods and for the recent development of so-called tensor methods mostly
studied in the machine learning community. Tensor methods
are based on the fact that in many latent variables models, low order moments contain enough information
to make the inverse problem efficiently solvable, with a solution that uses simple computational steps;
see \cite {AK2014} and references therein. \\

\noindent
The main goal of this paper is to extend such moment  methods in the setting of mixtures of binary regression models.
Let $X\in \mathbb{R}^{d}$ be the vector of covariates and $Y\in \lbrace 0,1\rbrace$ be the binary output.
A binary regression model assumes that for some link function $g$, the probability that $Y=1$ conditionally to $X=x$ is given by $g(\langle \beta , x \rangle +b)$, where $\beta\in  \mathbb{R}^{d}$ is the vector of regression coefficients and $b\in \mathbb{R}$ is the intercept. Popular examples of link functions are the logit link function where for any real $z$,  $g(z)=e^z/(1+e^z)$ and the probit link function where $g(z)=\Phi(z),$  with $\Phi$  the cumulative distribution function of the standard normal ${\cal N}(0,1)$. 
If now we want to modelise heterogeneous populations, let $K$ be the number of populations and $\omega=(\omega_1,\cdots,\omega_K)$ their weights such that $\omega_{j}\geq 0$, $j=1,\ldots,K$ and $\sum_{j=1}^{K}\omega_{j}=1$. Define, for $j=1,\ldots,K$, the regression coefficients in the $j$-th population by $\beta_{j}\in \mathbb{R}^{d}$ and the intercept in the $j$-th population by $b_{j}\in \mathbb{R}$. Let $\omega =(\omega_{1},\ldots,\omega_{K})$,   $b=(b_1,\cdots,b_K)$, $\beta=[\beta_{1} \vert \cdots,\vert \beta_K]$ the $d\times K$ matrix of regression coefficients and denote $\theta=(\omega,\beta,b)$.
The model of population mixture of binary regressions is given by:
\begin{equation}
\label{mixturemodel1}
 P_{\theta}(Y=1\vert X=x)=\sum^{K}_{k=1}\omega_k g(<\beta_k,x>+b_k).
\end{equation}

\noindent
Provable guarantees of estimation methods for mixtures of binary regressions are few, see Chapters 9 and 12 of \cite{MR3889980} and references therein, or \cite{AK2014b} and references therein. In  particular in \cite{AK2014b} it is proved that  cross moments up to order $3$ between the covariates and the output allow to recover the directions $\pm\beta_{j}/\|\beta_{j}\|$, $j=1,\ldots,K$.\\

\noindent
Our first main result is that cross moments up to order $3$  between the output and the regression variables are enough to recover all the parameters of the model. This holds for the probit link function, see Theorem \ref{thm1} and for general link functions under weak assumptions, see Theorem \ref{thm2}. These moment identifiability theorems are detailed in Section \ref{Identifiability}. \\

\noindent
Our second contribution is a new least squares moment method to estimate the parameter $\theta$ with theoretical guarantees.
Since empirical moments do not match exactly those defined with a distribution $P_{\theta}$ for some $\theta$,
the estimator $\widehat{\theta}$ is the one that minimizes the square distance of theoretical moments to empirical moments.
The estimation method is detailed in Section \ref{Estimation} together with the algorithm to compute the estimator.
We developed the associated R-package \verb"morpheus" available on the CRAN (\cite{Loum_Auder}).
Consistency and asymptotic normality of our least squares moment estimator are proved in Theorem \ref{thm3}.
We then compare experimentally our method to the maximum likelihood estimator (computed using the R-package
flexmix proposed in \cite{bg-papers:Gruen+Leisch:2007a}). The experiments are presented in Section \ref{Simu},
where we give indications about the advantages or disadvantages of both methods.\\

\noindent
We end the paper by investigating the identifiability in other population mixture models of binary regressions. 
For generalized linear models with outputs that are not binary, there is a large literature including estimation methods and practical algorithms such as in
\cite{Dang2017}, \cite{MR2523852}, \cite{Hennig2000}, \cite{INGRASSIA2016},
\cite{Ingrassia2015}, \cite{MARUOTTI2017}, \cite{Mazza2017}, \cite{JSSv086i02}, \cite{Punzo2017}, \cite{PUNZO2018};
see also the recent book \cite{MR3889980} and references therein.
However, the identifiability problem has not been fully explored in these articles and theoretical guarantees
are largely missing. In \cite{MR1108557}, the identifiability is proved for  finite mixtures of logistic
regression models where only the intercept varies with the population. In \cite{MR2476114}, finite mixtures
of multinomial logit models with varying and fixed effects are investigated, the proofs of the
identifiability results use the explicit form of the logit function. In \cite{MR3244553}, further non
parametric identifiability of the link function is proved, but only for models where the base exponential
models are identifiable for mixtures, which does not apply to binary outcome.
We provide in Section \ref{others} several identifiability results which are useful as a first step to
obtain theoretical guarantees in applications, such as in \cite{MR3086415}.
We prove that with a known  smooth enough link function, the directions of the regression vectors may be
recovered under the only assumption that they are distinct; see Theorem \ref{theobis1}.
Then, under the strengthened assumption that they are linearly independent,
we prove that the link function may be recovered in a non-parametrical way ; see Theorem \ref{theobis2}.
We then study the simultaneous use of continuous and categorical covariates and further propose assumptions
under which the parameters and the  link function may be recovered; see Theorem \ref{theobis3}.
We finally prove that, with longitudinal data having at least $3$ repetitions for each individual,
the whole model is identifiable under the weakest assumption that the regression directions are distinct;
see Theorem \ref{theobis4}.\\

\noindent
All the proofs of the different Theorems are postponed in Section \ref{proofs}.

%% file: identification.tex
\section{Moment identifiability}
\label{Identifiability}
\label{Notation}

We now define the cross moments that will be used to identify model (\ref{mixturemodel1}).\\
Let us denote $[n]$ the set $\lbrace 1,2,\ldots,n\rbrace$ and $e_i\in \mathbb{R}^d,$ the $i$-th canonical basis vector of $\mathbb{R}^d.$
Denote also $I_d\in \mathbb{R}^{d\times d}$ the identity matrix in $\mathbb{R}^{d}$.
The tensor product of $p$ euclidean spaces $\mathbb{R}^{d_i},\,\,i\in [p]$ is noted $\bigotimes_{i=1}^p\mathbb{R}^{d_i}.$
$T$ is called a real p-th order tensor if $T\in \bigotimes_{i=1}^p\mathbb{R}^{d_i}.$
For $p=1,$ $T$ is a vector in $\mathbb{R}^d$ and for $p=2$, $T$ is a $d\times d$ real matrix.
The $(i_1,i_2,\ldots,i_p)$-th coordinate of $T$ with respect the canonical basis is denoted $T[i_1,i_2,\ldots,i_p]$, $ i_1,i_2,\ldots,i_p\in [d].$
Define the cross moments between the response $Y$ and the covariable $X$:

\begin{itemize}
\item[--] $M_1(\theta):= E_{\theta}[Y.X],$
\item[--] $M_2(\theta):= E_{\theta}\Big[Y.\big(X\otimes X-\sum_{j\in[d]}Y.e_j\otimes e_j\big)\Big],$ and
\item[--] $M_3(\theta):= E_{\theta}\Big[Y\big(X\otimes X\otimes X-\sum_{j\in[d]}\big[X\otimes e_j\otimes e_j+e_j\otimes X\otimes e_j+e_j\otimes e_j\otimes X\big]\big)\Big]$.
\end{itemize}

\noindent
We assume that the random variable $X$ has a Gaussian distribution.
We now focus on the situation where $X\sim \mathcal{N}(0,I_d)$, $I_d$ being the identity $d\times d$ matrix.
All results may be easily extended to the situation where $X\sim \mathcal{N}(m,\Sigma)$, $m\in  \mathbb{R}^{d}$,
$\Sigma$ a positive and symetric $d\times d$ matrix by writing $X=m+\sigma^{1/2}\tilde{X}$ with $\tilde{X}\sim \mathcal{N}(0,I_d)$,
and estimating empirically $m$ and $\Sigma$.
We shall explain below how the results can be extended to the case where $X$ has any smooth enough distribution with support $\mathbb{R}^d$.\\

\noindent
To prove our moment identifiability result, we shall use the following assumptions:
\begin{itemize}
\item[--] (H1)
The vectors $\beta_{1},\ldots,\beta_{K}$ are linearly independent and the weights are positive: $\omega_{k}>0$, $k=1,\ldots,K$.
\item[--] (H2)
The link function $g$ is strictly increasing from $0$ at $-\infty$ to $1$ at $+\infty$, it has continuous derivatives till order $4$,  strictly decreasing first derivative on $[0,+\infty[$, and it satisfies
$$\forall z\in  \mathbb{R},\;g(z)+g(-z)=1.$$
\item[--]
(H3) There exists a neighborhood $\cal O$ of $(0,0)$ in $\mathbb{R} _+^\star\times  \mathbb{R}$
and functions $L_s$, $s=1,2,3,$  such that for $s=1,2,3,$
$\int_{ R} L_s(z) e^{-z^{2}/2} dz < +\infty$ , and
$\forall$ $z\in  \mathbb{R},$ $\forall$ $(\lambda,b)\in {\cal O}$, we have 
$$(|z|+1)\left|\frac{\partial g^{(s+1)}}{\partial \lambda}(\lambda z+b)\right|\leq L_s(z).$$ 
\end{itemize}
\noindent
Notice that (H1) implies that $d\geq K$, and that (H2) and (H3) hold in particular for the logistic link function and the probit link function.

\begin{theorem}[Probit identifiability]
\label{thm1}
If (H1) holds and if $g$ is the probit link function, one may recover $K$ and $\theta=(\omega,\beta,b)$ from the knowledge of $M_{1}(\theta)$, $M_{2}(\theta)$ and $M_{3}(\theta)$.
\end{theorem}
For general link functions, identifiability holds at least in an open set.
The following identifiability result and Theorem \ref{thm1} are proved in  Section \ref{proof:theo1}.
\begin{theorem}[General identifiability]
\label{thm2}
If (H1), (H2), (H3) hold and $g^{(3)}(0)\neq 0,$  there exist $L>0$ and $B>0$ such that as soon as  $\|\beta_{k}\| <L$ and $|b_{k}|<B$  for all $k=1,\ldots,K$,
then one may recover $K$ and $\theta=(\omega,\beta,b)$ from the knowledge of $M_{1}(\theta)$, $M_{2}(\theta)$ and $M_{3}(\theta)$.
\end{theorem}
Since the proof uses Taylor expansions, it only proves the existence of {\it small enough} positive $L$ and $B$ such that the result holds. When the link function $g$ is the logit function,
we investigated numerically the function defined in (\ref{esp}) in the proof for which the property of being one-to-one is enough to deduce identifiability,
and we found  that Theorem \ref{thm2} seems to be true at least with $L=8$ and $B=8$.\\

\noindent
For population mixture models of linear regressions with continuous outcomes, it is only needed that the vectors of regression coefficients are distinct.
However, for population mixture models of binary regressions, identifiability results exist only for the logit link function, see \cite{MR1108557} and \cite{MR2476114}.
In Section \ref{others} we prove that the directions of the regression vectors are identifiable under the only assumption that they are distinct. One consequence of Theorem \ref{thm2} is that it provides provable guarantees for the maximum likelihood estimator so that in particular it is consistent.\\

\noindent
The Gaussian assumption for the distribution of $X$ may be relaxed by defining score-adapted cross-moments in the same way as in Section 4 of \cite{AK2014b}.


\section{The least squares moment estimator}
\label{Estimation}

\subsection{Definition of the estimator}
\label{estimator}
In Section \ref{Identifiability} we showed that the parameters can be recovered by matching the cross-moments till order $3$.
Those moments are unknown, so we estimate them empirically using:
\begin{eqnarray*}
  \widehat{M}_1&=&\frac{1}{n}\sum_{i=1}^n Y_i X_i\\
  \widehat{M}_2&=&\frac{1}{n}\sum_{i=1}^n\Big[Y_i (X_i\otimes X_i-\sum_{j\in[d]} e_j\otimes e_j )\Big]\\
  \widehat{M}_3&=&\frac{1}{n}\sum_{i=1}^n \Big[Y_i (X_i\otimes X_i\otimes X_i-
    \sum_{j\in[d]}\big[X_i\otimes e_j\otimes e_j+e_j\otimes X_i\otimes e_j+e_j\otimes e_j\otimes X_i\big] )\Big].
\end{eqnarray*}
It is not possible to match the empirical moments exactly, so we use a least-squares estimator. Let $\Theta$ be the set of parameters, and define for all $\theta$:
$$
Q_n(\theta)=\sum_{j\in[d]} \Big\lbrace\widehat{M}_1[j]-M_1(\theta)[j]\Big\rbrace^2+\sum_{j,k\in [d]}\Big\lbrace\widehat{M}_2[j,k]-M_2(\theta)[j,k]\Big\rbrace^2
+\sum_{j,k,l\in [d]}\Big\lbrace\widehat{M}_3[j,k,l]-M_3(\theta)[j,k,l]\Big\rbrace^2.
$$
The least squares moment method (LSMM) defines the estimator as follows:
\begin{eqnarray}
\label{eq9}
\widehat{\theta}_n=\underset{\theta\in \Theta}{\text{argmin}}\,\,Q_n(\theta).
\end{eqnarray}
\noindent
The terms in the sum to compute $Q_n(\theta)$ contribute unevenly to the result, because
there are more indices combinations for higher order moments, and thus more (sub-)terms.
Given this observation, pondering each term might lead to improvements in the estimation.
Following \cite{Hansen1982} we propose the following generalized least squares moment method (GLSMM).
Define for all $\theta$ and $i=1,\ldots,n$:
$$\tilde M_i(\theta) = Y_i
\begin{pmatrix}
  X_i\\
  X_i\otimes X_i-\sum_{j\in[d]} e_j\otimes e_j\\
  X_i\otimes X_i\otimes X_i-\sum_{j\in[d]}\big[X_i\otimes e_j\otimes e_j+e_j\otimes X_i\otimes e_j+e_j\otimes e_j\otimes X_i\big]\\
\end{pmatrix}
-
\begin{pmatrix}
  M_1(\theta)\\
  M_2(\theta)\\
  M_3(\theta)\\
\end{pmatrix}
$$
as a column vector, by flattening the tensorial products (for example put column $j$ between indices $j \times d$ and $(j+1) \times d - 1$).
Let $W$ be a symmetric and positive definite square matrix of size $d + d^2 + d^3$, 
set
$$Q^{W}_n(\theta)=\left(\frac{1}{n} \sum_{i=1}^{n} {}^t \negthinspace \tilde M_i(\theta) \right) W \left(\frac{1}{n} \sum_{i=1}^{n}\tilde M_i(\theta)\right), $$
and define
\begin{eqnarray}
\label{eq9bis}
\widehat{\theta}^{W}_n=\underset{\theta\in \Theta}{\text{argmin}}\,\,Q^{W}_n(\theta).
\end{eqnarray}
Notice that when $W$ is the identity matrix$I_d$, $\widehat{\theta}^{I_d}_n=\widehat{\theta}_n$.

\subsection{Asymptotic properties}
\label{asymptotics}
Define $\theta^{\star}$ as the true unknown parameter.\\
To prove the consistency of the estimator, identifiability of the model is obviously necessary.
Thus, if the link function is not probit, then we choose $\Theta$ such that all parameters in $\Theta$ satisfy  $\|\beta_{k}\| <L$ and $|b_{k}|<B$
for all $k=1,\ldots,K$, $L$ and $B$ being defined in Theorem \ref{thm2}.\\
To prove the limiting Gaussian distribution of our least-squares moment estimator, we shall need more assumptions.
For $j=1,\ldots,5$, let $G_{j}$ be the $K\times K$ diagonal matrix having the $ E[g^{(j)}\left( \langle \beta_{k}^{\star}, X\rangle + b_{k}^{\star}\right)]$'s on the diagonal. 

\begin{itemize}
  \item[--] (H4) All diagonal coefficients of $G_{3}$ are non zero.
  \item[--] (H5) All diagonal coefficients of $G_{1}G_{3}-G_{2}^{2}$ are non zero.
\end{itemize}

Let $q=K(2+d)-1$ be the dimension of $\theta$.
Let $Z_{n}$ be the gradient vector of $Q_n(\theta)$ at $\theta = \theta^*$.
$Z_{n}$ is a linear combination of empirical moments, see the explicit formula in Section \ref{proof:theo3}.
Define $\Gamma(\theta^{\star})$ the $q\times q$ matrix which is the variance matrix of $\sqrt{n}Z_{n}$.
For each $\theta$, define $\tilde M(\theta)$ as the expectation  of $\tilde M_i(\theta)$, that is
$$\tilde M(\theta) =
\begin{pmatrix}
  M_1(\theta^{\star})\\
  M_2(\theta^{\star})\\
  M_3(\theta^{\star})\\
\end{pmatrix}-
\begin{pmatrix}
  M_1(\theta)\\
  M_2(\theta)\\
  M_3(\theta)\\
\end{pmatrix}
$$
as a column vector in the same way as $\tilde M_i (\theta)$. 
Define also the $q \times q$ matrix $V(\theta)$ such that for all $r_{1},r_{2}=1,\ldots,q$:
$$
V_{r_1r_2}(\theta)= 2  \frac{{}^t \negthinspace\partial\tilde M(\theta)}{\partial  \theta_{r_1}}W \frac{\partial \tilde M(\theta)}{\partial \partial \theta_{r_2}}.
$$
\begin{theorem}
\label{thm3}
Assume that (H1), (H2), (H3) hold,  that $\Theta$ is compact and that $\theta^{\star}\in\Theta$. Then $\widehat{\theta}^{W}_n$ is consistent.\\
If moreover (H4) and (H5) hold, then $V(\theta^{\star})$ is invertible, and $\sqrt{n}\left(\widehat{\theta}^{W}_n-\theta^\star\right)$ converges
in distribution under $P_{\theta^{\star}}$ to a centered Gaussian distribution with variance $\Sigma(\theta^{\star})=V(\theta^{\star})^{-1}\Gamma(\theta^{\star}) V(\theta^{\star})^{-1}$.
\end{theorem}
The proof of Theorem \ref{thm3} is detailed in Section \ref{proof:theo3} and follows the usual analysis of the asymptotic behavior of $Z$-estimators,
the more delicate part of the proof being to prove that the limiting value  $V(\theta^{\star})$ of the Hessian of $Q_{n}(\theta)$ is invertible.\\
To apply this Theorem one may estimate consistently $\Sigma(\theta^{\star})$ by plug-in.\\

\noindent
A natural question is the choice of the matrix $W$. It is proved in \cite{Hansen1982} that the matrix $W$ minimizing the asymptotic variance of the estimator is given by $W(\theta^{\star})$
$W(\theta) = ( \E [ \tilde M_i(\theta) {}^t \negthinspace \tilde M_i(\theta) ])^{-1}$.
This optimal matrix can be estimated empirically as follows:
$$\hat W(\hat \theta) = \left( \frac{1}{n} \sum_{i=1}^{n} \tilde M_i(\hat \theta) {}^t \negthinspace \tilde M_i(\hat \theta) \right)^{-1} \, ,$$
with $\hat \theta$ the parameters estimated from the LSMM for example, or another GLSMM (with another choice of the matrix $W$).
$W(\theta)$ can also be re-estimated as many times as needed in a loop algorithm.

\subsection{Algorithm}
\label{algorithm}
The estimator $\widehat{\theta}_{n}$ is computed by using the representation of the regression vectors $\beta_{k}$
through their direction $\mu_{k}$ and norm $\lambda_{k}$, $k=1,\ldots,K$.
In a first step, we compute a preliminary estimate of $[\mu_{1},\ldots,\mu_{K}]$ using a spectral method.
In a second step, we minimize $Q_{n}^{W_{\texttt{init}}}$ using usual optimization methods with $W_{\texttt{init}} = Id$, or any user-defined matrix.
$W$ can then be recomputed using the optimized parameters, and a final minimization of $Q_{n}^W$ follows.
It would be possible to iterate this procedure again, but in practice it did not lead to improved results.
The preliminary estimator for the directions is used as initial point for the directions in the optimization procedure.

\begin{table}[H]
\centering
\begin{tabular}{p{13.5cm}}
\hline
Algorithm $M3LS$: Estimation of all parameters\\
\hline
\textbf{input}: $X,\,Y, \,K,\,g$ \\
$1:$ Estimate the directions $\mu_{1},\ldots,\mu_{K}$ using Algorithm $InitDir$ \\
$2:$ Optimize $Q_n^{W_{\texttt{init}}} (\theta)$ using the estimators of $1.$ as initial directions. Stop here if $W_{\texttt{init}}$ is considered good enough \\
$3:$ Re-compute $W$ using the optimized parameters $p$, $\beta$ and $b$ \\
$4:$ Execute step $2.$ again \\
\textbf{Output}: The estimated parameter $\widehat{\theta}$ \\
\hline
\end{tabular}
\label{tab01}
\end{table}

The preliminary estimation of the directions is based on the spectral method. For any vector $z\in \mathbb{R}^{p}$, define $B(z)$ the $d\times d$ matrix such that
$$B(z)[i,j]:=\sum_{s=1}^{d}M_{3}(\theta)[i,j,s]z_{s},$$
so that, using Lemma \ref{lem_rewrite}, we get
$$B(z)=\sum_{k=1}^{K} \omega_k\lambda_k^3\mathbb{E}[g^{(3)}\big(\lambda_k<X,\mu_k>+b_k\big)] \langle \mu_k,z\rangle  \mu_k^{\otimes 2}.$$
It is proved in \cite{Afsari2008} that it is possible to recover the directions by joint diagonalisation of $B(z_{1}), \ldots B(z_{P})$ for distinct vectors $z_{1},\ldots,z_{P}$, $P\geq 2$.
Joint diagonalisation of $B(z_{1}), \ldots B(z_{P})$ means finding a matrix $V$ such that the matrices $VB(z_p)V^T$ are the most diagonal possible.
The normalized vectors $\mu_{1},\ldots,\mu_{K}$ are obtained up to sign and label switching by taking the first $K $ vectors of $V^{-1}.$
Let us denote $U$ the matrix of these $K$ vectors. Let $O=U_{\ast}^{-1}M_1(\theta)\in \mathbb{R}^{K},$  with  $U_{\ast}^{-1}$ the general inverse of $U.$
The real numbers $\omega_k\lambda_k\mathbb{E}[g'\big(\lambda_k\langle X,\mu_k\rangle+b_k\big)]$, $k=1,\ldots,K$, are given up to sign by the elements of $O.$
Since they 
are positive, the sign of the $\mu_k$'s are obtained by multiplying $-1$ all the vectors associated to the negative values of $O.$\\
In pratice, the vectors $\mu_{1},\ldots,\mu_{K}$ are estimated using the joint diagonalisation method applied to the matrices $\widehat{B}(z_{p})$, $p=1,\ldots,P$, computed using $\widehat{M}_{3}$.

\begin{table}[H]
\centering
\begin{tabular}{p{13.5cm}}
\hline
Algorithm $InitDir$: Joint diagonalisation algorithm to estimate the directions\\
\hline
\textbf{input}: $X,\,Y, K$ \\
$1$: Estimate the cross moments $\widehat{M}_1,$ $\widehat{M}_2$ and $\widehat{M}_3$ as explained in section   \ref{asymptotics}\\
$2$: Choose vectors $\lbrace z_1,z_2\ldots,z_P\rbrace\subseteq  \mathbb{R}d$ (for instance: the canonical basis $ e_1,e_2,\ldots,e_P$ of $ \mathbb{R}^d$ if $P=d$)\\
$3$: Compute $\widehat{B}(z_{p})$ for all $p\in\lbrace 1,2,\ldots,P\rbrace$\\
$4$: Joint diagonalisation: compute $V$ such that   $V\widehat{B}(z_{p})V^T$ are the most diagonal possible \\
$5$: Compute $U=V^{-1}[1:K]$ the K-first vectors of $V^{-1}$ (by ordering the diagonal values in decreasing absolute value)\\
$6$: Compute $O=ginv(U) \widehat{M}_1 $\\
$7$: Multiply by $-1$ all the vectors of U corresponding to the negative values of $O$ $U[,O<0]=-U[,O<0]$\\
\textbf{Output}: The preliminar estimators of $\mu_{1},\ldots,\mu_{K}$ \\
\hline
\end{tabular}
\label{tab02}
\end{table}

%% file: simulations.tex
\section{Simulations}
\label{Simu}

We first present the developed R-package implementing the generalized least-squares moment method (GLSMM) in Section \ref{pkg}.
Then, in Section \ref{compareperfs}, the algorithms GLSMM and the maximum likelihood method (MLM)
are compared in terms of performance (accuracy and computation time).
We also compare in Section \ref{initialization} the initialization as indicated in the introduction below,
with an initialization based on random starts.

\subsection{\textsf{R} \textbf{package}}
\label{pkg}

The developed R-package is called \verb"morpheus" \cite{Loum_Auder} and is divided into two main parts:
\begin{enumerate}
  \item the computation of the matrix $\mu$ (containing the normalized columns of $\beta$),
    based on the empirical cross-moments as described in section \ref{algorithm};
  \item the optimization of all parameters (including $\beta$), using the initially estimated
    directions as a starting point.
\end{enumerate}
The former is a straightforward translation of the mathematical formulas (file R/computeMu.R),
while the latter calls R constrOptim() method on the objective function expression and its
derivative (file R/optimParams.R). For usage examples, please refer to the package help.

\subsection{\textbf{Comparison of the GLSMM and MLM algorithms}}
\label{compareperfs}

In this section, we compare experimentally our GLSMM algorithm (morpheus package \cite{Loum_Auder}) to the
MLM algorithm (with flexmix package \cite{bg-papers:Gruen+Leisch:2007a} which is a reference for this kind of
estimation, using an iterative algorithm to maximize the log-likelihood).
We present at first estimation errors and we compare computational time for the two different approaches.
The parameters for the simulations are chosen arbitrarily as indicated below;
these should be discovered by the algorithms (GLSMM, and the MLM algorithm).\\

Experiment 1 (2 dimensions):
\begin{align*}
	K &= 2\\
	p &= (0.5, 0.5)\\
	b &= (-0.2, 0.5)\\
	\beta &=
		\begin{pmatrix}
		1 & 3\\
		-2 & 1
		\end{pmatrix}
\end{align*}

Experiment 2 (5 dimensions):
\begin{align*}
	K &= 2\\
	p &= (0.5, 0.5)\\
	b &= (-0.2, 0.5)\\
	\beta &=
		\begin{pmatrix}
			1 & 2\\
			2 & -3\\
			-1 & 0\\
			0 & 1\\
			3 & 0
		\end{pmatrix}
\end{align*}

Experiment 3 (10 dimensions):
\begin{align*}
	K &= 3\\
	p &= (0.3, 0.3, 0.4)\\
	b &= (-0.2, 0, 0.5)\\
	\beta &=
		\begin{pmatrix}
			1 & 2 & -1\\
			2 & -3 & 1\\
			-1 & 0 & 3\\
			0 & 1 & -1\\
			3 & 0 & 0\\
			4 & -1 & 0\\
			-1 & -4 & 2\\
			-3 & 3 & 0\\
			0 & 2 & 1\\
			2 & 0 & -2
		\end{pmatrix}
\end{align*}

An experiment consists in a data generation step, followed by the parameters estimation.
The same link function is used for these two stages.
For all three experiments we use both logit and probit links, which amounts to a total of six different computations.
Within each experiment, the sample size varies from $5\times 10^3$ to $10^6$ ($5\times10^3,$
$10^4,$ $10^5,$ $5\times10^5$ and $10^6$ ).
Computations are always run on the same data both for GLSMM and MLM.
Results are aggregated over $N=100$ replications.\\

\noindent \textbf{Estimation errors}.

To compare the estimation errors between GLSMM and MLM, we present in Tables \ref{logit_res} and
\ref{probit_res} the Mean Absolute Error summed over all components of each parameter.
That is to say, for $p$ (resp. $b$) the error showed is $\sum_{k=1}^{K} |p_k - \hat p_k|$ (resp. $\sum_{k=1}^{K} |b_k - \hat b_k|$).
Concerning $\beta$, we show the error column per column: $\texttt{err}_k = \sum_{i=1}^{d} |\beta_{i,k} - \hat \beta_{i,k}|$ with $k \in [1,K]$.
The parameters are estimated by averaging the outputs of $N=100$ replications for different values of $n$ in each Experiment.
Then, for each method we remove the 2\% quantile of extreme (largest L1 norm) values, which are rather rare and would bias the visual summary too much.\\

Performances are similar, with an advantage to MLM when using the logit link function,
and to GLSMM when using the probit link function.
Since MLM is asymptotically optimal, this is a good result.
Note that the first line for the MLM algorithm ends with very high values (third experiment, $d=10$):
this is due to a high percentage of failed runs (about one third), so that the 2\% quantile matches less than two runs.\\

\begin{table}[H]
 \begin{center}
  \begin{tabular}{ccccc|cccc|ccccc}
   \hline \\
   & \multicolumn{4}{c}{Experiment 1} &  \multicolumn{4}{c}{Experiment 2} & \multicolumn{5}{c}{Experiment 3}\\
   \cline{2-14} \\
   $n$ & $p$ & $\beta_1$ & $\beta_2$ & $b$ & $p$ & $\beta_1$ & $\beta_2$ & $b$ & $p$ & $\beta_1$ & $\beta_2$ & $\beta_3$ & $b$\\
   \hline \\
   \multirow{2}{*}{$5.10^3$} & 1.3e-2 & 1.1e+0 & 1.3e+0 & 2.6e-1 & 1.8e-2 & 2.1e+0 & 1.6e+0 & 3.9e-1 & 4.9e-2 & 5.2e+0 & 3.8e+0 & 1.2e+1 & 7.8e+0\\
	 & 4.0e-3 & 7.7e-2 & 3.1e-2 & 8.6e-3 & 6.0e-4 & 2.7e-1 & 7.0e-2 & 6.3e-3 & 1.3e-1 & 1.6e+2 & 1.8e+2 & 1.4e+0 & 4.3e+0\\
	 \\
	 \multirow{2}{*}{$10^4$} & 1.2e-2 & 3.5e-1 & 4.1e-1 & 7.0e-2 & 2.5e-2 & 7.1e-1 & 9.7e-1 & 4.6e-1 & 3.3e-2 & 3.9e+0 & 3.4e+0 & 7.4e+0 & 2.7e+0\\
	 & 2.0e-3 & 3.1e-2 & 7.1e-2 & 1.4e-2 & 1.8e-3 & 3.3e-2 & 4.2e-2 & 6.9e-3 & 1.3e-1 & 2.2e+0 & 1.8e+0 & 4.8e-1 & 7.3e-2\\
	 \\
	 \multirow{2}{*}{$10^5$} & 7.6e-3 & 3.6e-2 & 5.1e-2 & 9.7e-3 & 5.8e-2 & 4.1e-1 & 2.6e-1 & 3.5e-2 & 1.7e-2 & 1.1e+0 & 6.8e-1 & 2.0e+0 & 2.7e-1\\
	 & 7.7e-4 & 6.6e-3 & 6.7e-3 & 3.2e-3 & 5.9e-5 & 2.2e-2 & 1.7e-2 & 2.4e-3 & 1.3e-1 & 4.6e-2 & 9.6e-2 & 9.3e-2 & 5.5e-3\\
	 \\
	 \multirow{2}{*}{$5.10^5$} & 3.6e-3 & 2.7e-2 & 4.2e-2 & 1.1e-2 & 1.9e-2 & 2.0e-1 & 8.0e-2 & 9.4e-3 & 1.6e-2 & 9.1e-1 & 1.0e+0 & 7.9e-1 & 7.6e-2\\
	 & 1.2e-4 & 3.4e-3 & 8.1e-3 & 5.0e-3 & 2.3e-4 & 5.0e-3 & 6.1e-3 & 3.9e-3 & 1.3e-1 & 2.8e-2 & 1.7e-2 & 1.7e-2 & 2.5e-3\\
	 \\
	 \multirow{2}{*}{$10^6$} & 6.3e-4 & 8.9e-3 & 1.3e-2 & 4.4e-3 & 7.0e-2 & 5.3e-1 & 4.0e-1 & 1.9e-2 & 7.5e-3 & 7.8e-1 & 8.7e-1 & 2.7e-1 & 7.8e-2\\
	 & 2.8e-5 & 6.1e-4 & 2.3e-3 & 1.5e-3 & 7.0e-5 & 2.5e-3 & 7.0e-3 & 2.5e-3 & 1.3e-1 & 5.4e-2 & 2.0e-2 & 1.0e-2 & 3.4e-3\\ \\
	 \hline
  \end{tabular}
 \end{center}
 \caption{\textit{logit} link function. Summed errors for GLSMM (top number) and MLM (bottom number), for increasing values of $n$.
  Parameters estimations are averaged over $N=100$ runs.}
 \label{logit_res}
\end{table}

\begin{table}[H]
 \begin{center}
  \begin{tabular}{ccccc|cccc|ccccc}
   \hline \\
   & \multicolumn{4}{c}{Experiment 1} &  \multicolumn{4}{c}{Experiment 2} & \multicolumn{5}{c}{Experiment 3}\\
   \cline{2-14} \\
   $n$ & $p$ & $\beta_1$ & $\beta_2$ & $b$ & $p$ & $\beta_1$ & $\beta_2$ & $b$ & $p$ & $\beta_1$ & $\beta_2$ & $\beta_3$ & $b$\\
   \hline \\
	 \multirow{2}{*}{$5.10^3$} & 3.0e-2 & 2.0e+0 & 1.1e+0 & 2.6e-1 & 5.4e-3 & 3.0e+0 & 2.5e+0 & 5.4e-1 & 5.2e-2 & 1.6e+0 & 1.3e+0 & 8.2e+0 & 6.0e+0\\
	 & 3.3e-3 & 1.8e-2 & 4.9e-2 & 5.0e-3 & 2.0e-3 & 6.8e-2 & 1.1e-1 & 2.2e-2 & 1.4e-1 & 1.0e+1 & 7.4e+1 & 1.5e+0 & 2.0e-1\\
   \\
	 \multirow{2}{*}{$10^4$} & 1.1e-2 & 1.0e+0 & 1.6e+0 & 2.2e-1 & 2.4e-2 & 1.6e+0 & 2.1e+0 & 4.2e-1 & 2.2e-2 & 2.1e+0 & 2.0e+0 & 7.4e+0 & 2.4e+0\\
	 & 1.2e-3 & 8.2e-2 & 1.7e-1 & 4.0e-2 & 4.5e-4 & 1.7e-1 & 1.0e-1 & 1.1e-2 & 1.4e-1 & 4.0e+0 & 2.2e+0 & 2.5e-1 & 5.1e-2\\
   \\
	 \multirow{2}{*}{$10^5$} & 6.3e-3 & 1.3e-1 & 6.2e-2 & 1.1e-2 & 2.9e-2 & 5.5e-1 & 2.8e-1 & 3.5e-2 & 1.9e-2 & 2.1e+0 & 1.7e+0 & 3.0e+0 & 3.7e-1\\
	 & 1.2e-3 & 5.1e-1 & 1.0e+0 & 2.1e-1 & 5.4e-5 & 8.4e-1 & 7.8e-1 & 9.6e-2 & 1.3e-1 & 1.3e+1 & 1.2e+1 & 8.0e+0 & 5.4e-1\\
   \\
	 \multirow{2}{*}{$5.10^5$} & 3.7e-3 & 3.8e-2 & 8.5e-2 & 9.1e-3 & 3.1e-2 & 1.6e-1 & 2.4e-1 & 1.0e-2 & 1.2e-2 & 2.3e+0 & 1.9e+0 & 2.7e+0 & 1.9e-1\\
	 & 2.4e-4 & 1.0e+0 & 2.0e+0 & 4.1e-1 & 5.4e-5 & 1.1e+0 & 1.1e+0 & 1.3e-1 & 1.3e-1 & 1.7e+1 & 1.6e+1 & 1.1e+1 & 7.4e-1\\
	 \\
	 \multirow{2}{*}{$10^6$} & 4.9e-3 & 4.8e-2 & 1.2e-1 & 1.5e-2 & 2.0e-2 & 1.5e-1 & 3.4e-2 & 9.1e-3 & 1.0e-2 & 2.6e+0 & 1.9e+0 & 1.6e+0 & 1.3e-1\\
 	 & 1.2e-4 & 1.3e+0 & 2.6e+0 & 5.2e-1 & 1.0e-5 & 1.6e+0 & 1.5e+0 & 1.9e-1 & 1.3e-1 & 1.7e+1 & 1.6e+1 & 1.1e+1 & 7.4e-1\\ \\
   \hline
  \end{tabular}
 \end{center}
 \caption{\textit{probit} link function. Summed errors for GLSMM (top number) and MLM (bottom number), for increasing values of $n$.
  Parameters estimations are averaged over $N=100$ runs.}
 \label{probit_res}
\end{table}

\noindent \textbf{Computational time}.

The running time of LSMM (GLSMM without the generalized least squares step) does not depend directly on $n$,
since after computing the empirical moments it operates on matrices of size at most $O(d^3 \times K)$.
So, we observed excellent running times in this case, outperforming MLM by far especially for $n=10^6$, as seen on Table \ref{timing_nooptim}.
The times where averaged over 100 runs using random true parameters, with the logit link function.
Only one core is available for each run (they are executed in parallel).
The top numbers correspond to the GLSMM algorithm, and the numbers below correspond to the flexmix package (MLM).\\

\begin{table}[H]
 \begin{center}
  \begin{tabular}{c|c|c|c}
   \hline \\
   & $n=10^4$ & $n=10^5$ & $n=10^6$\\
    \cline{2-4} \\
    \multirow{2}{*}{$d=2$} & 0.77 & 0.77 & 1.1 \\
    & 5.2 & 36.8 & 288.3 \\
   \\
	 \multirow{2}{*}{$d=5$} & 1.3 & 1.5 & 2.4 \\
    & 5.8 & 49.4 & 501.4 \\
   \\
	 \multirow{2}{*}{$d=10$} & 5.6 & 6.8 & 8.7 \\
    & 7.2 & 76.3 & 756.2 \\ \\
    \hline
  \end{tabular}
 \end{center}
  \caption{Average running time over 100 runs \emph{with a fixed $W$ matrix} on one core, for the logit link, for different values of $n$ and $d$. Top numbers: GLSMM; bottom numbers: MLM.}
 \label{timing_nooptim}
\end{table}

However, the accuracy of the initial method in high dimension was clearly bad, so the final step computing and using the empirical matrix $W$ is required.
This final step turns out to be quite costly, so the timings in this case get closer, with still an advantage to GLSMM as seen on Table \ref{timing_optim}.
This time we used four cores per run to benefit from a parallel loop in the computation of $W$, in addition to the parallelization of the runs.
As we can see this also helps the flexmix package.\\

\begin{table}[H]
 \begin{center}
  \begin{tabular}{c|c|c|c}
   \hline \\
   & $n=10^4$ & $n=10^5$ & $n=10^6$\\
    \cline{2-4} \\
    \multirow{2}{*}{$d=2$} & 16.2 & 13.3 & 9.8 (*) \\
    & 2.3 & 16.1 & 116.6 \\
   \\
	 \multirow{2}{*}{$d=5$} & 9.9 & 13.0 & 24.5 \\
    & 2.6 & 23.5 & 181.3 \\
   \\
	 \multirow{2}{*}{$d=10$} & 26.3 & 52.4 & 220.3 \\
    & 3.2 & 30.2 & 298.3 \\ \\
   \hline
  \end{tabular}
 \end{center}
 \caption{Average running time over 100 runs \emph{with optimisation of $W$} on four cores, for the logit link, for different values of $n$ and $d$. Top numbers: GLSMM; bottom numbers: MLM.}
 \label{timing_optim}
\end{table}

\noindent (*) This anomaly could be explained by a quicker convergence when $n$ increases,
doing more than compensating for the slightly longer computation time for the moments and the matrix $W$.
When $d$ is larger, the optimisation step is no longer the only factor.\\

We can imagine to benefit from the initial speed of the algorithm by estimating $W$ on a reduced dataset,
before the parameters are estimated from the full dataset starting from a custom $W_{\texttt{init}}$.
Overall, the results are quite promising, especially if we find a way to optimize the computation of $W$.


\subsection{\textbf{Random initialization}}
\label{initialization}

Instead of choosing $\mu$ as a starting point for the optimization stage,
we could just sample a few points at random in the unit sphere $S^{d-1}$ and keep
the initial point leading to the smallest value of $Q_n(\theta)$.
So we compare here these two settings, choosing arbitrarily three random points on the sphere.
This version of GLSMM is written GLSMM3 in this section.
Table \ref{randomStarts} summarizes the results, showing in each cell the total error for $p$, $b$ or a column of the matrix $\beta$.
First two lines correspond to the logit link, respectively for GLSMM and GLSMM3, whereas last two lines correspond to the probit link function.
This time, no extreme values were removed because no serious outliers were observed.\\

Experimental setup:
\begin{itemize}
  \item $n=10^5$ sample points,
  \item parameters according to the experiments 1, 2 and 3 respectively ($d=2$, then $5$ and $10$).
\end{itemize}

\begin{table}[H]
 \begin{center}
  \begin{tabular}{ccccc|cccc|ccccc}
   \hline \\
   & \multicolumn{4}{c}{Experiment 1} &  \multicolumn{4}{c}{Experiment 2} & \multicolumn{5}{c}{Experiment 3}\\
   \cline{2-14} \\
   link & $p$ & $\beta_1$ & $\beta_2$ & $b$ & $p$ & $\beta_1$ & $\beta_2$ & $b$ & $p$ & $\beta_1$ & $\beta_2$ & $\beta_3$ & $b$\\
   \hline \\
	 \multirow{2}{*}{logit} & 6.2e-3 & 2.8e-2 & 8.5e-3 & 5.7e-3 & 1.2e-2 & 4.6e-2 & 1.0e-1 & 4.2e-2 & 1.6e-2 & 7.3e-1 & 1.6e+0 & 2.7e+0 & 2.8e-1\\
   & 2.6e-2 & 5.7e-2 & 5.4e-2 & 3.9e-3 & 1.7e-3 & 7.4e-2 & 1.4e-1 & 4.1e-2 & 1.2e-2 & 5.8e-1 & 2.0e+0 & 1.5e+0 & 2.3e-1\\
	 \\
	 \multirow{2}{*}{probit} & 5.6e-3 & 5.8e-2 & 1.1e-2 & 6.5e-3 & 1.5e-2 & 1.3e-1 & 4.4e-1 & 1.7e-1 & 1.1e-2 & 1.2e+0 & 5.6e-1 & 3.9e+0 & 3.6e-1\\
   & 2.6e-2 & 2.2e-1 & 1.8e-1 & 3.1e-2 & 6.0e-4 & 1.2e-1 & 6.1e-1 & 6.8e-2 & 8.7e-3 & 4.6e-1 & 3.2e-1 & 4.0e+0 & 3.8e-1\\ \\
	 \hline
  \end{tabular}
 \end{center}
 \caption{Summed errors for GLSMM (top numbers) and GLSMM3 (bottom numbers), for $n=10^5$ and both logit and probit link functions.
  Parameters estimations are averaged over $N=100$ runs.}
 \label{randomStarts}
\end{table}

The random starts version of GLSMM is slightly more accurate, although it is not easily spotted on the tables.
The total sums of errors for the logit link are respectively 5.5 (GLSMM) and 4.7 (GLSMM3), and
respectively 6.9 and 6.4 for the probit link function.
However, using three random starting points is costly, because of the current computational burden to obtain the weights matrix $W$.
Moreover, using only one random starting point (GSLMM1) leads to a total sum of errors quite higher than GSLMM:
9.4 for the logit link, and 9.1 for the probit link.
As a conclusion, using $\mu$ as a starting point for the optimization step is a good compromise.

%% file: extensions.tex
\section{Some other identifiability results}
\label{others}
In this section, we provide several further identifiability results for various population mixture models of binary regressions.\\

\noindent
Identifiability of a model is a first step to obtain theoretical guarantees for practical estimation procedures. Our identifiability results open the way to build estimators for which theoretical guarantees could be obtained. In particular, for parametric maximum likelihood estimators in mixture models for which algorithms already exist, consistency is a consequence of our identifiability theorems  by applying the usual theory.

\subsection{Continuous covariates}

We first consider the setting considered in the previous sections where
$$
E\left( Y \vert X\right)= \sum^{K}_{k=1}\omega_k g(\langle \beta_k,X \rangle+b_k).
$$
We assume that for all $k$, $\omega_k \geq 0$, that $\sum_{k=1}^{K}\omega_{k}=1$, and that  $g$ takes value in $(0,1)$. We show below that the directions of the regression vectors may be recovered as soon as they are distinct, even if the link function is unknown.\\
Denote $ P_{g,\omega,\beta,b}$ the probability distribution of  $(X,Y)$, with $\omega=(\omega_1,\cdots,\omega_K)$, $\beta=[\beta_1|,\cdots,|\beta_K] \in \mathbb{R}^{d\times K}$, and 
$b=(b_1,\cdots,b_K) \in \mathbb{R}^{K}$.
When $g$ is unknown obviously it is needed to fix origin and scale, we choose to fix $g(0)$ and $g(1)$ (with no loss of  generality). Denote $\mu_{k}=\beta_{k}/\|\beta_{k}\|$ and $\lambda_{k}=\|\beta_{k}\|$, $k=1,\ldots,K$, so that $\beta_{k}=\lambda_{k} \mu_{k}$.
\\
We introduce the assumptions:
\begin{itemize}
\item[--] (S1)
The support of the law of $X$ is $ \mathbb{R}^d$.
\item[--]  (S2)
For all $j\neq k$,  $\mu_{j}\neq \mu_{k}$ and $\mu_{j}\neq -\mu_{k}$.
\item[--]  (S3)
The function $g: \mathbb{R} \rightarrow ]0,1[$ is increasing, has limit $0$ in $-\infty$,  limit $1$ in $+\infty$, and it is continuously derivable with derivative having   limit $0$ in $-\infty$ and in $+\infty$. Also, $g(0)<g(1)$ are fixed.
\end{itemize}
\noindent
\underline{Remark}: There is no assumption on $K$ with respect to $d$.

\begin{theorem}
\label{theobis1}
Under assumptions (S1), (S2) and (S3), knowledge of $ P_{g,\omega,\beta,b}$  allows to recover $K$ and $\mu_{1},\ldots,\mu_{K}$.
\end{theorem}
\noindent
The proof of this theorem is given in Section \ref{sec:theobis1}.

\noindent
Under the more stringent assumption that the regression vectors are linearly independent, it is possible to recover all parameters and the link function.
\begin{itemize}
\item[--] (S2bis)
The vectors $\mu_{1},\ldots, \mu_{K}$ are linearly independent.
\end{itemize}
\noindent
\underline{Remark}: (H2bis) implies that $K\leq d$.

\begin{theorem}
\label{theobis2}
Under assumptions (S1), (S2bis) and (S3), the mixture model is identifiable: the knowledge of $ P_{g,\omega,\beta,b}$  allows to recover $K$, $g$, $\omega$, $\beta$ and $b$.
\end{theorem}
\noindent
The proof of this theorem is given in Section \ref{sec:theobis2}.

\subsection{Continuous and categorical covariates}

We now consider the situation where part of the covariates are catergorical, we denote them  $Z$, and  $\{z_{1},\ldots,z_{m}\}\subset \mathbb{R}^{d'}$ their possible values. We still denote $X\in \mathbb{R}^{d}$ the continuous covariates. Now
$$
E\left( Y \vert X, Z\right)= \sum^{K}_{k=1}\omega_k g(\langle \beta_k,X \rangle+ \langle \gamma_k,Z \rangle + b_k),
$$
and we denote $ P_{g,\omega,\beta,\gamma, b}$ the probability distribution of $(X,Y)$, with $\omega=(\omega_1,\cdots,\omega_K)$, $\beta=[\beta_1|,\cdots,|\beta_K] \in \mathbb{R}^{d\times K}$, $\gamma=[\gamma_1|,\cdots,|\gamma_K] \in \mathbb{R}^{d'\times K}$, and
$b=(b_1,\cdots,b_K) \in \mathbb{R}^{K}$.
We introduce
\begin{itemize}
\item[--]  (S4)
The matrix  $\left(\begin{array}{ll}1 & z_{1}^{T}\\
1 & z_{2}^{T}\\
\vdots & \vdots \\
1 & z_{m}^{T}
\end{array}\right)$ is full rank.
\end{itemize}

\noindent
\underline{Remark:}
(S4) implies that $d'+1 \leq m$.\\

\noindent
The continuous covariates allow to identify $g$.

\begin{theorem}
\label{theobis3}
Under assumptions (S1), (S2bis), (S3) and (S4), the model is identifiable: the knowledge of $ P_{g,\omega,\beta,\gamma,b}$  allows to recover $K$, $g$, $\omega$, $\beta$, $\gamma$ et $b$.
\end{theorem}
\noindent
The proof of this theorem is given in Section \ref{sec:theobis3}.

\subsection{Longitudinal observations}

We now consider the situation where for each individual $Y$, conditional to the membership of a population, we have $p$ independent experiments with several covariates $X_{1},\ldots,X_{p}$. Thus the random variable $Y$ has dimension $m$, and
$$
E\left( Y \vert X, Z\right)= \sum^{K}_{k=1}\omega_k
\left(g(\langle \beta_k,X_{j} \rangle+ \langle \gamma_k,Z_{j} \rangle + b_k)\right)_{1\leq j \leq p}.
$$

\noindent
As soon as the number of experiments is at least $3$, we do not need the linear independence of the regression vectors to get identifiability. 

\begin{theorem}
\label{theobis4}
Assume that $p\geq 3$.
If (S1), (S2), (S3) and (S4) hold, then  the model is identifiable: the knowledge of $ P_{g,\omega,\beta,\gamma,b}$  allows to recover $K$, $g$, $\omega$, $\beta$, $\gamma$ and $b$.
\end{theorem}
\noindent
The proof of this theorem is given in Section \ref{sec:theobis4}.

%% file: proofs.tex
\section{Proofs}
\label{proofs}

\subsection{Proof of Theorem \ref{thm1} and Theorem \ref{thm2}}
\label{proof:theo1}

\noindent
Let, for $k=1,\ldots,K$, $\lambda_k = \|\beta_k \| $ and $\mu_{k}=\beta_{k}/\|\beta_k \| $.
Using Stein's identity, Anandkumar et al. (\cite{AK2014b}) prove the following lemma: 
\begin{lemma}[\cite{AK2014b}]
\label{lem_rewrite}
Under (H3)  the moments can be rewritten:
\begin{eqnarray*}
		M_1(\theta)&=&\sum_{k=1}^K\omega_k\lambda_k E[g'\big(\lambda_k\langle X,\mu_k\rangle+b_k\big)]\;\mu_k ,
		\\
		M_2(\theta)&=&\sum_{k=1}^K\omega_k\lambda_k^2 E[g''\big(\lambda_k<X,\mu_k>+b_k\big)]\;\mu_k\otimes \mu_k ,
		\\
		M_3(\theta)&=&\sum_{k=1}^K\omega_k\lambda_k^3 E[g^{(3)}\big(\lambda_k<X,\mu_k>+b_k\big)]\;\mu_k\otimes \mu_k\otimes \mu_k .
\end{eqnarray*}
\end{lemma}
Under (H1), we see by Lemma \ref{lem_rewrite} that $K$ is the rank of $M_{2}(\theta)$. \\
It is proved in \cite{AK2014b} that we can recover the $\mu_{k}$'s up to sign from the knowledge of $ M_{2}(\theta)$ and $ M_{3}(\theta)$, but since  $M_{1}(\theta)$ is a linear combination of the $\mu_{k}$'s with positive coefficients, under (H1) the knowledge of  $M_{1}(\theta)$ allows to recover the signs. It is then seen that using $M_{1}(\theta)$, $M_{2}(\theta)$ and $M_{3}(\theta)$, one may recover the $3$-uples
$$
\left(\omega_{k}E[ g'\left( \langle \beta_{k}, X\rangle + b_{k}\right)]\lambda_{k};\omega_{k} E[g"\left( \langle \beta_{k}, X\rangle + b_{k}\right)]\lambda_{k}^{2};\omega_{k} E[g^{(3)}\left( \langle \beta_{k}, X\rangle + b_{k}\right)]\lambda_{k}^{3} \right),
$$
$k=1,\ldots,K.$ Thus, one gets identifiability as soon as the function from $]0,+\infty[\times   \mathbb{R} \times ]0,+\infty[$ to its image that associates $(\omega, b,\lambda)$ to
$$
\left(\omega\lambda \int g'(\lambda z+b)e^{-z^{2}/2}dz; \omega\lambda^{2} \int g"(\lambda z+b)e^{-z^{2}/2}dz;\omega \lambda^{3}\int g^{(3)}(\lambda z+b)e^{-z^{2}/2}dz\right)
$$
is one-to-one. Using integration by parts this is equivalent to the fact that
the function from $]0,+\infty[\times  \mathbb{R}\times ]0,+\infty[$ to its image that associates $(\omega,b,\lambda)$ to
$$
\lambda \left(\omega \int g'(\lambda z+b)e^{-z^{2}/2}dz; \omega \int z g'(\lambda z+b)e^{-z^{2}/2}dz;\omega \int z^{2} g'(\lambda z+b)e^{-z^{2}/2}dz\right)
$$
is one-to-one. This is again equivalent to the fact that the function from $\mathbb{R}\times ]0,+\infty[ $ to its image that associates $(b,\lambda)$ to
$$
\left( \frac{\int z g'(\lambda z+b)e^{-z^{2}/2}dz}{\int g'(\lambda z+b)e^{-z^{2}/2}dz};\frac{ \int z^{2} g'(\lambda z+b)e^{-z^{2}/2}dz}{\int g'(\lambda z+b)e^{-z^{2}/2}dz}\right)
$$
is one-to-one. For any $(b,\lambda)\in \mathbb{R}\times ]0,+\infty[$, define 
\begin{eqnarray}
\label{loi}
dQ_{(b,\lambda)}(z)=\frac{g'(\lambda z+b)e^{-z^{2}/2}}{\int g'(\lambda z+b)e^{-z^{2}/2}dz}dz.
\end{eqnarray}
Then it is equivalent to prove that the knowledge of 
\begin{eqnarray}
\label{esp}
\left(E_{(b,\lambda)}(Z);E_{(b,\lambda)}(Z^{2})\right):=\left(\int z dQ_{(b,\lambda)}(z); \int z^{2} dQ_{(b,\lambda)}(z)\right)
\end{eqnarray}
implies the knowledge of $(b,\lambda)$. 

\subsubsection{End of the proof of Theorem \ref{thm1}}
When the link function $g$ is probit, then $g'(z)=\frac{1}{\sqrt{2\pi}}e^{-z^2/2}$. Replacing in equation (\ref{loi}), we have
$$
dQ_{(b,\lambda)}(z)=\frac{e^{-\frac{1}{2}\left((\lambda z+b)^2+z^{2}\right)}}{\int e^{-\frac{1}{2}\left((\lambda z+b)^2+z^{2}\right)}dz}dz,
$$
which after some computations leads to
 $$Q_{(b,\lambda)}={\cal N}\left(-\frac{\lambda b}{\lambda^{2}+1};\frac{1}{\lambda^{2}+1}\right).$$ 
Its first two moments are then given by
$$(\alpha_1, \alpha_2 ) :=\left( E_{(b,\lambda)}(Z); E_{(b,\lambda)}(Z^{2})\right)=\left(-\frac{\lambda b}{\lambda^{2}+1};\frac{\lambda^2 b^2+\lambda^2+1}{(\lambda^{2}+1)^2}\right).$$
We can then recover $b$ and $\lambda$ by 
$$b=-\alpha_1\frac{(\lambda^2+1)}{\lambda}$$ and $$
\lambda=\sqrt{(\alpha_2-\alpha_1^2)^{-1}-1}.$$ 

\subsubsection{End of the proof of Theorem \ref{thm2}} 
It remains to prove that for some $B>$ and $L>0$, if $(H2)$ holds, then, the function
that associates $(b,\lambda)\in ]-B,B[\times ]0,L[$ to
$
\left(E_{(b,\lambda)}(Z);E_{(b,\lambda)}(Z^{2})\right)$ is one-to-one on its image.\\
Using (H3) and integration by parts, we get that for $(b,\lambda)$ in a neighborhood of $(0,0)$
\begin{eqnarray}
\label{eq_bij}
\left( E_{(b,\lambda)}(Z); E_{(b,\lambda)}(Z^{2})\right)=\left(\frac{\lambda\int g''(\lambda z+b)e^{-z^{2}/2}dz}{\int g'(\lambda z+b)e^{-z^{2}/2}dz};1+\frac{\lambda^2\int  g^{(3)}(\lambda z+b)e^{-z^{2}/2}dz}{\int g'(\lambda z+b)e^{-z^{2}/2}dz}\right).
\end{eqnarray}
From (H2), one gets that
\begin{itemize}
\item
(P1) The function $g'$ is positive and satisfies $g'(x)=g'(-x)$ for all $x\in \mathbb{R}$,
\item
(P2) The function $g"$ satisfies  $g"(x)=-g"(-x)$ for all $x\in \mathbb{R}$  and $g"(x)< 0$ for $x>0$,
\item
(P3) $g'(0)>0$ ,
\item
(P4) $g''(0)=g^{(4)}(0)=0$.
\end{itemize}
Let us define the functions $K_s,$ $s=1,2,3$ such that
$$\begin{array}{ccccl}
 K_s & : &  \mathbb{R}\times  ]0,+\infty [ & \to &  \mathbb{R} \\
 & & (b,\lambda) & \mapsto & K_s(b,\lambda) =\int g^{(s)}(\lambda z+b)e^{-z^2/2}dz
\end{array}$$
Using (H3), the functions $K_s,$ $s=1,2,3$, are differentiable in a neighborhood of $(0,0)$ and Taylor expansion  writes:
\begin{eqnarray}
\label{dt1}
K_s(b,\lambda)=K_s(0,0) + \langle \nabla K_s(0,0), (b,\lambda)\rangle+o(\lambda^2+b^2).
\end{eqnarray}
Now
\begin{eqnarray}
\label{dp3}
\frac{\partial K_s}{\partial \lambda}(0,0)=\int zg^{(s+1)}(0)e^{-z^2/2}dz
\end{eqnarray}
and
\begin{eqnarray}
\label{dp4}
\frac{\partial K_s}{\partial b}(0,0)=\int g^{(s+1)}(0)e^{-z^2/2}dz
\end{eqnarray}
so that
\begin{eqnarray}
\label{dt2}
K_s(b,\lambda)=g^{(s)}(0)\int e^{-z^2/2}dz +g^{(s+1)}(0)\int (\lambda z+b) e^{-z^2/2}dz +o(\lambda^2+b^2).
\end{eqnarray}
Using $(P4)$ and (\ref{dt2}), we have
\begin{equation}
\label{der1}
\int g'(\lambda z+b)e^{-z^2/2}dz=\sqrt{2\pi} g'(0)+o(\lambda^2+b^2),
\end{equation}
\begin{eqnarray}
\label{der2}
\int g''(\lambda z+b)e^{-z^2/2}dz=\sqrt{2\pi} g^{(3)}(0)b+o(\lambda^2+b^2),
\end{eqnarray}
and
\begin{eqnarray}
\label{der3}
\int g^{(3)}(\lambda z+b)e^{-z^2/2}dz=\sqrt{2\pi} g^{(3)}(0)+o(\lambda^2+b^2).
\end{eqnarray}
Therefore, replacing (\ref{der1}) to (\ref{der3}) in (\ref{eq_bij}), we get
$$
E_{(b,\lambda)}(Z)=\frac{g^{(3)}(0)}{g'(0)}\lambda b + o(\lambda^{2}+b^{2})
$$
and
$$
E_{(b,\lambda)}(Z^{2})=1 + \frac{g^{(3)}(0)}{g'(0)}\lambda^{2} + o(\lambda^{2}+b^{2}),
$$
which easily leads to
$$
\lambda^{2}= \frac{g'(0)}{g^{(3)}(0)}\left(E_{(b,\lambda)}(Z)^{2}-1\right)+ o(\vert E_{(b,\lambda)}(Z)^{2}-1 \vert + \vert E_{(b,\lambda)}(Z)\vert )
$$
and
$$
\lambda b=\frac{g'(0)}{g^{(3)}(0)}E_{(b,\lambda)}(Z) + o(\vert E_{(b,\lambda)}(Z)^{2}-1 \vert + \vert E_{(b,\lambda)}(Z)\vert ).
$$
This proves that the function $(b,\lambda)\mapsto (E_{(b,\lambda)}(Z),E_{(b,\lambda)}(Z^2))$ is invertible in a neighborhood of $(0,0)$.

\subsection{Proof of Theorem \ref{thm3}}
\label{proof:theo3}
By the law of large numbers, $Q^{W}_{n}(\theta)$ converges to 
$$
Q^{W}(\theta):=  {}^t \negthinspace \tilde M (\theta)  W \tilde M(\theta).
$$
Define
$$S_n=\sup _{\theta \in \Theta}\Big\vert Q_n(\theta)-Q(\theta)\Big\vert.$$
Since $Q(\theta)$ has $\theta^{\ast}$ as unique minimum (up to label switching), to prove the consistency of $\widehat{\theta}_{n}$, it is enough to prove that $S$ converges to $0$ in probability, see Theorem 5.7 in \cite{vaart}.
Let $\lambda_{max}$ be the largest eigenvalue of $W$ and$D=(d+d^{2}+d^{3})^{2}$.
We easily get
\begin{eqnarray*}
S_n &\leq&\displaystyle \lambda_{max}D\sum_{j\in[d]} \Big(\Big\vert\hat{M}_1[j]-M_1(\theta^\star)[j]\Big\vert\Big)\Big( \big\vert\hat{M}_1[j]\big\vert+\big\vert M_1(\theta^{\ast})[j]\big\vert+2 \sup _{\theta \in \Theta}\big\vert M_1(\theta)[j]\big\vert\Big)\\ 
&+&\displaystyle\lambda_{max}D\sum_{j,k\in[d]} \Big(\Big\vert\hat{M}_2[j,k]-M_2(\theta^\star)[j,k]\Big\vert\Big)\Big( \big\vert\hat{M}_2[j,k]\big\vert+\big\vert M_2(\theta^{\ast})[j,k]\big\vert+2 \sup _{\theta \in \Theta}\big\vert M_2(\theta)[j,k]\big\vert\Big)\\ 
&+&\displaystyle\lambda_{max}D\sum_{j,k,l\in[d]} \Big(\Big\vert\hat{M}_3[j,k,l]-M_3(\theta^\star)[j,k,l]\Big\vert\Big)\Big( \big\vert\hat{M}_3[j,k,l]\big\vert+ \big\vert M_3(\theta^{\ast})[j,k,l]\big\vert+2\sup _{\theta \in \Theta}\big\vert M_3(\theta)[j,k,l]\big\vert\Big),
\end{eqnarray*}
and since the functions $\theta\mapsto\,\, M_r(\theta),\,\,r=1,2,3$ are continuous and $\Theta$ is compact, then there exist $c_1,\,\,c_2$ and $c_3$ such that
\begin{eqnarray*}
S_n &\leq&\displaystyle\lambda_{max}D\sum_{j\in[d]} \Big(c_1+\big\vert\hat{M}_1[j]\big\vert\Big)\Big(\Big\vert\hat{M}_1[j]-M_1(\theta^\star)[j]\Big\vert\Big)
\\&+&\displaystyle\lambda_{max}D\sum_{j,k\in[d]}\Big(c_2+\big\vert\hat{M}_2[j,k]\big\vert\Big) \Big(\Big\vert\hat{M}_2[j,k]-M_2(\theta^\star)[j,k]\Big\vert\Big)\\ 
&+&\displaystyle\lambda_{max}D\sum_{j,k,l\in[d]} \Big(c_3+\big\vert\hat{M}_3[j,k,l]\big\vert\Big)\Big(\Big\vert\hat{M}_3[j,k,l]-M_3(\theta^\star)[j,k,l]\Big\vert\Big)
\end{eqnarray*}
which converges to $0$ by the law of large numbers, which ends the proof of the consistency of $\widehat{\theta}_{n}$.\\

Let us define $Z_n$ as
$
Z_n(\theta)=\nabla_{\theta}Q_n(\theta)$.
The $r-th$ coordinate of $Z_n(\theta)$ can be obtained by
$$
\frac{\partial Q_n(\theta)}{\partial \theta_r} =-2  {}^t \negthinspace\frac{\partial\tilde M(\theta)}{\partial \theta_r}W  \left(\frac{1}{n} \sum_{i=1}^{n}\tilde M_i(\theta)\right).
$$
Using Taylor expansion, we get
\begin{eqnarray}
Z_n(\hat{\theta}_n)=Z_n(\theta^\star)+\int_0^1D_1Z_n\big[\theta^{\star}+t(\hat{\theta}_n-\theta^\star)\big]\Big(\hat{\theta}_n-\theta^{\star}\Big)dt
\end{eqnarray}
where $D_1Z_n$ is the first derivative matrix of $Z_n$. Since $Z_n(\hat{\theta}_n)=0,$ we have
\begin{eqnarray}
\label{eq1}
-\sqrt{n}Z_n(\theta^\star)=\left[\int_0^1D_1Z_n\big[\theta^{\star}+t(\hat{\theta}_n-\theta^\star)\big]dt\,\, \right]\sqrt{n}\Big(\hat{\theta}_n-\theta^\star\Big)
\end{eqnarray}
Applying the central limit theorem and  the delta method we get that $\sqrt{n}Z_n(\theta^\star)$ is asymptotically Gaussian.

The $(r_1,r2)-th$ coordinate of $D_1Z_n(\theta)=\nabla^2_{\theta}Q_n(\theta)$ are given by 
$$
\frac{\partial^2 Q_n(\theta)}{\partial \theta_{r_1}\partial \theta_{r_2}}=-2  {}^t \negthinspace\frac{\partial^{2}\tilde M(\theta)}{\partial  \theta_{r_1}\partial \theta_{r_2}}W  \left(\frac{1}{n} \sum_{i=1}^{n}\tilde M_i(\theta)\right)+ V_{r_1r_2}(\theta)
$$
%
with
$$
V_{r_1r_2}(\theta)=2  {}^t \negthinspace\frac{\partial\tilde M(\theta)}{\partial  \theta_{r_1}}W \frac{\partial \tilde M(\theta)}{\partial \partial \theta_{r_2}}.
$$
It is not difficult to prove that 
$\int_0^1D_1Z_n\big[\theta^{\star}+t(\hat{\theta}_n-\theta^{\star})\big]dt$ converges in probability to $V(\theta^{\star})$ 
 so that the proof is completed by showing that the matrix $V=V(\theta^\star)$ is invertible. 

$V$ is a $q\times q$ matrix with $q=K(2+d)-1$. Let $U\in  \mathbb{R}^{q}$. We shall denote the coordinates of $U$ according to the parameters.
Since $W$ is a symmetric and positive definite square matrix , using the form of $V$ we get that $U^T V U=0$ if and only if:
\begin{eqnarray}
\label{ft}
U^{T} D M_{1}(\theta^{\star})[j]=0,\;j=1,\ldots,q,
\end{eqnarray}
and
\begin{eqnarray}
\label{st}
U^{T} D M_{2}(\theta^{\star})[j,l]=0,\;j,l=1,\ldots,q,
\end{eqnarray}
and
\begin{eqnarray}
\label{tt}
U^{T} D M_{3}(\theta^{\star})[j,l,m]=0,\;j,l,m=1,\ldots,q.
\end{eqnarray}
Here, $D M_{\cdot}[]$ is the gradient vector of the involved coordinate of $M_{\cdot}$. Denote $U(\beta_{k})$ the $d$-dimensional vector involving the coordinates of $U$ according to parameter $\beta_{k}$. Denote $\overline{0}$ the $d$-dimensional zero vector, $\overline{0}\otimes \overline{0}$ the  $d\times d$-dimensional zero matrix and $\overline{0}\otimes \overline{0}\otimes \overline{0}$ the  $d\times d\times d$-dimensional zero third order tensor. Then, the equation (\ref{ft}) can be rewritten as:
\begin{equation}
\label{ft2}
U^TDM_1(\theta^{\star})[j]=\sum_{k=1}^{K-1}U(\omega_k)\frac{\partial M_1(\theta^{\star})[j]}{\partial \omega_k}+ \sum_{k=1}^{K}U(b_k)\frac{\partial M_1(\theta^{\star})[j]}{\partial b_k}
+ \sum_{k=1}^{K}\sum_{m=1}^dU(\beta_{mk})\frac{\partial M_1(\theta^{\star})[j]}{\partial \beta_{mk}},
\end{equation}
the equation (\ref{st}) can be rewritten as
\begin{equation}
\label{st2}
U^TDM_2(\theta^{\star})[j,l]=\sum_{k=1}^{K-1}U(\omega_k)\frac{\partial M_2(\theta^{\star})[j,l]}{\partial \omega_k}+ \sum_{k=1}^{K}U(b_k)\frac{\partial M_2(\theta^{\star})[j,l]}{\partial b_k}+ \sum_{k=1}^{K}\sum_{m=1}^dU(\beta_{mk})\frac{\partial M_2(\theta^{\star})[j,l]}{\partial \beta_{mk}},
\end{equation}
and the equation (\ref{tt}) can be rewritten as
\begin{equation}
\label{tt2}
U^TDM_3(\theta^{\star})[j,l,m]=\sum_{k=1}^{K-1}U(\omega_k)\frac{\partial M_3(\theta^{\star})[j,l,m]}{\partial \omega_k}+ \sum_{k=1}^{K}U(b_k)\frac{\partial M_3(\theta^{\star})[j,l,m]}{\partial b_k}
+ \sum_{k=1}^{K}\sum_{s=1}^dU(\beta_{sk})\frac{\partial M_3(\theta^{\star})[j,l,m]}{\partial \beta_{sk}}.
\end{equation}
Using the fact that $\sum_{k=1}^d\omega_k=1,$ the first terms of the equations (\ref{ft2}) to (\ref{tt2}) are rewritten as:
\begin{eqnarray}
\sum_{k=1}^{K-1}U(\omega_k)\frac{\partial M_1(\theta^{\star})[j]}{\partial \omega_k}&=&\sum_{k=1}^{K-1}U(\omega_k)\Big\lbrace  E\left[g'\left(\langle X,\beta^{\star}_k\rangle+b^{\star}_k\right)\right].\beta^{\star}_k(j)\nonumber\\
\label{ft3}
&-& E\left[g'\left(\langle X,\beta^{\star}_K\rangle+b^{\star}_K\right)\right].\beta^{\star}_K(j)\Big\rbrace,
\end{eqnarray}
\begin{eqnarray}
\sum_{k=1}^{K-1}U(\omega_k)\frac{\partial M_2(\theta^{\star})[j,l]}{\partial \omega_k}&=&\sum_{k=1}^{K-1}U(\omega_k)\Big\lbrace  E\left[g''\left(\langle X,\beta^{\star}_k\rangle+b^{\star}_k\right)\right].\beta^{\star}_k(j)\beta_k(l)\nonumber\\
\label{st3}
&-& E\left[g''\left(\langle X,\beta^{\star}_K\rangle+b^{\star}_K\right)\right].\beta^{\star}_K(j)\beta^{\star}_K(l)\Big\rbrace
\end{eqnarray}
and
\begin{eqnarray}
\sum_{k=1}^{K-1}U(\omega_k)\frac{\partial M_3(\theta^{\star})[j,l,m]}{\partial \omega_k}&=&\sum_{k=1}^{K-1}U(\omega_k)\Big\lbrace  E\left[g^{(3)}\left(\langle X,\beta^{\star}_k\rangle+b^{\star}_k\right)\right].\beta^{\star}_k(j)\beta^{\star}_k(l)\beta^{\star}_k(m)\nonumber\\
\label{tt3}
&-& E\left[g^{(3)}\left(\langle X,\beta^{\star}_K\rangle+b^{\star}_K\right)\right].\beta^{\star}_K(j)\beta^{\star}_K(l)\beta^{\star}_K(m)\Big\rbrace
\end{eqnarray}
respectively. Likewise the seconds terms of equations (\ref{ft2}) to (\ref{tt2}) are rewritten as:
\begin{eqnarray}
\label{ft4}
\sum_{k=1}^{K}U(b_k)\frac{\partial M_1(\theta^{\star})[j]}{\partial b_k}=\sum_{k=1}^K\omega^{\star}_kU(b_k) E\left[g''\left(\langle X,\beta^{\star}_k\rangle+b^{\star}_k\right)\right].\beta_{k}^{\star}(j),
\end{eqnarray}
\begin{eqnarray}
\label{st4}
\sum_{k=1}^{K}U(b_k)\frac{\partial M_2(\theta^{\star})[j,l]}{\partial b_k}=\sum_{k=1}^K\omega^{\star}_kU(b_k) E\left[g^{(3)}\left(\langle X,\beta^{\star}_k\rangle+b^{\star}_k\right)\right].\beta^{\star}_{k}(j)\beta^{\star}_{k}(l)
\end{eqnarray}
and 
\begin{eqnarray}
\label{tt4}
\sum_{k=1}^{K}U(b_k)\frac{\partial M_3(\theta^{\star})[j,l,m]}{\partial b_k}=\sum_{k=1}^K\omega^{\star}_kU(b_k) E\left[g^{(4)}\left(\langle X,\beta^{\star}_k\rangle+b^{\star}_k\right)\right].\beta^{\star}_k(j)\beta^{\star}_k(l)\beta^{\star}_k(m)
\end{eqnarray}
respectively. Derivating with respect to the $\beta_k$'s coordinates and using  Stein's identity, the last terms of equations (\ref{ft2}) to (\ref{tt2}) are rewritten as: 

\begin{eqnarray}
\sum_{k=1}^{K}\sum_{m=1}^dU(\beta_{mk})\frac{\partial M_1(\theta^{\star})[j]}{\partial \beta_{mk}}&=&\sum_{k=1}^K\omega^{\star}_k E\left[g^{(3)}\left(\langle X,\beta^{\star}_k\rangle+b^{\star}_k\right)\right]\langle \beta^{\star}_k,U(b_k)\rangle\beta^{\star}_k(j)\nonumber\\
\label{ft5}
&+& \sum_{k=1}^K\omega^{\star}_k E\left[g'\left(\langle X,\beta^{\star}_k\rangle+b^{\star}_k\right)\right]U(\beta_k(j)),
\end{eqnarray}

\begin{eqnarray}
\sum_{k=1}^{K}\sum_{m=1}^dU(\beta_{mk})\frac{\partial M_2(\theta^{\star})[j,l]}{\partial \beta_{mk}}&=&\sum_{k=1}^K\omega^{\star}_k E\left[g^{(4)}\left(\langle X,\beta^{\star}_k\rangle+b^{\star}_k\right)\right]\langle \beta^{\star}_k,U(b_k)\rangle\beta^{\star}_k(j)\beta_k(l)\nonumber\\
&+& \sum_{k=1}^K\omega^{\star}_k E\left[g''\left(\langle X,\beta^{\star}_k\rangle+b^{\star}_k\right)\right]\Big\lbrace \beta^{\star}_k(j)U(\beta_k(l))\nonumber\\
\label{st5}
&+&\beta^{\star}_k(l)U(\beta_k(j))\Big\rbrace,
\end{eqnarray}
and

\begin{eqnarray}
\sum_{k=1}^{K}\sum_{s=1}^dU(\beta_{sk})\frac{\partial M_3(\theta^{\star})[j,l,m]}{\partial \beta_{sk}}&=&\sum_{k=1}^K\omega^{\star}_k E\left[g^{(5)}\left(\langle X,\beta^{\star}_k\rangle+b^{\star}_k\right)\right]\langle \beta^{\star}_k,U(b_k)\rangle\beta^{\star}_k(j)\beta_k(l)\beta^{\star}_k(s)\nonumber\\
&+& \sum_{k=1}^K\omega^{\star}_k E\left[g^{(3)}\left(\langle X,\beta^{\star}_k\rangle+b^{\star}_k\right)\right]\Big\lbrace \beta^{\star}_k(j)\beta^{\star}_k(l)U(\beta^{\star}_k(s))\nonumber\\
\label{tt5}
&+&\beta^{\star}_k(j)U(\beta_k(l))\beta^{\star}_k(s)+U(\beta_k(j))\beta^{\star}_k(l)\beta^{\star}_k(s)\Big\rbrace
\end{eqnarray}
respectively. Then using equations (\ref{ft2}) to (\ref{tt5}), we can rewrite equation (\ref{ft}) as:
\begin{eqnarray}
\overline{0}&=&\sum_{k=1}^{K-1} U(\omega_{k})\Big\lbrace E\big[ g'\left( \langle X,\beta^{\star}_{k}\rangle + b^{\star}_{k}\right)\big]\beta^{\star}_{k}- E\big[ g'\left( \langle X,\beta^{\star}_{K}\rangle + b^{\star}_{K}\right)\big]\beta_{K}\Big\rbrace\nonumber\\
&+&\sum_{k=1}^{K}\omega^{\star}_{k}U(b_{k}) E\big[ g''\left( \langle X,\beta^{\star}_{k}\rangle + b^{\star}_{k}\right)\big]\beta^{\star}_{k}+\sum_{k=1}^{K}\omega_{k} E\big[ g^{(3)}\left( \langle X,\beta^{\star}_{k}\rangle + b^{\star}_{k}\right)\big]\langle U(\beta_{k}),\beta^{\star}_{k}\rangle \beta^{\star}_{k}\nonumber\\
&+&\sum_{k=1}^{K}\omega^{\star}_{k} E\big[ g'\left( \langle X,\beta^{\star}_{k}\rangle + b^{\star}_{k}\right)\big]U(\beta_{k}),
\label{1}
\end{eqnarray}
rewrite equation (\ref{st}) as:
\begin{eqnarray}
\overline{0}\otimes \overline{0}
&=&\sum_{k=1}^{K-1} U(\omega_{k})\Big[ E\big[ g''\left( \langle X,\beta^{\star}_{k}\rangle + b^{\star}_{k}\right)\big]\beta^{\star}_{k}\otimes \beta_{k}- E\big[ g''\left( \langle X,\beta^{\star}_{K}\rangle + b^{\star}_{K}\right)\big]\beta^{\star}_{K}\otimes \beta_{K}\Big]\nonumber\\
&+&\sum_{k=1}^{K}\omega_{k}U(b_{k}) E\big[ g^{(3)}\left( \langle X,\beta^{\star}_{k}\rangle + b^{\star}_{k}\right)\big]\beta^{\star}_{k}\otimes \beta^{\star}_{k}\nonumber\\&+&\sum_{k=1}^{K}\omega^{\star}_{k} E\big[ g^{(4)}\left( \langle X,\beta^{\star}_{k}\rangle + b^{\star}_{k}\right)\big]\langle U(\beta_{k}),\beta^{\star}_{k}\rangle \beta^{\star}_{k}\otimes \beta^{\star}_{k} \nonumber\\
&+&\sum_{k=1}^{K}\omega^{\star}_{k} E\big[ g''\left( \langle X,\beta^{\star}_{k}\rangle + b^{\star}_{k}\right)\big]\Big(U(\beta_{k})\otimes \beta^{\star}_{k}+\beta^{\star}_{k}\otimes U(\beta_{k})\Big),
\label{2}
\end{eqnarray}
and rewrite equation (\ref{tt}) as:
\begin{eqnarray}
\overline{0}^{\otimes 3}
  &=&\sum_{k=1}^{K-1} U(\omega_{k})\Big[ E\big[ g^{(3)}\left( \langle X,\beta^{\star}_{k}\rangle + b^{\star}_{k}\right)\big]{\beta^{\star}_{k}}^{\otimes 3}- E\big[ g^{(3)}\left( \langle X,\beta^{\star}_{K}\rangle + b^{\star}_{K}\right)\big]{\beta^{\star}_{K}}^{\otimes 3}\Big]\label{3}\\
  &+&\sum_{k=1}^{K}\omega^{\star}_{k}U(b_{k}) E\big[ g^{(4)}\left( \langle X,\beta^{\star}_{k}\rangle + b^{\star}_{k}\right)\big]{\beta^{\star}_{k}}^{\otimes 3}+\sum_{k=1}^{K}\omega^{\star}_{k} E\big[ g^{(5)}\left( \langle X,\beta^{\star}_{k}\rangle + b^{\star}_{k}\right)\big]\langle U(\beta_{k}),\beta^{\star}_{k}\rangle {\beta^{\star}_{k}}^{\otimes 3} \nonumber\\
&+&\sum_{k=1}^{K}\omega_{k} E\big[ g^{(3)}\left( \langle X,\beta^{\star}_{k}\rangle + b^{\star}_{k}\right)\big]\Big(U(\beta_{k})\otimes \beta^{\star}_{k}\otimes \beta_{k}+\beta^{\star}_{k}\otimes U(\beta_{k})\otimes \beta^{\star}_{k}+ \beta^{\star}_{k}\otimes \beta^{\star}_{k}\otimes U(\beta_{k})\Big).
\nonumber
\end{eqnarray}
We shall first prove that the vectors $U(\beta_{1}),\ldots,U(\beta_{K})$ all belong to the linear space spanned by  $\beta_{1},\ldots,\beta_{K}$. 

Let $W$ be any vector that is orthogonal to this  linear space. By multiplying   (\ref{3}) on the right by $W$, and by using the fact that  $\beta_{1},\ldots,\beta_{K}$ are linearly independent by (H1), we get that
$$
\forall k=1,\ldots,K,\;\omega^{\star}_{k} (G_{3})_{k}\langle U(\beta_{k}),W\rangle=0.
$$
Using (H1) we have $\omega^{\star}_{k}>0$, $k=1,\ldots,K$ so that we get 
$$
\forall k=1,\ldots,K,\;  (G_{3})_{k}\langle U(\beta_{k}),W\rangle=0.
$$
Then, if (H4) holds, we get that for any $k$ and any $W$, $\langle U(\beta_{k}),W\rangle =0$, which proves that  the vectors $U(\beta_{1}),\ldots,U(\beta_{K})$ all belong to the linear space spanned by  $\beta_{1},\ldots,\beta_{K}$.
Let $B$ be the $d\times K$ matrix having $\beta_{1},\ldots,\beta_{K}$ as colomn vectors. Let $U(\beta)$ be the $d\times K$ matrix having $U(\beta_{1}),\ldots,U(\beta_{K})$ as colomn vectors. We thus have that there exists a  $K \times K$ matrix $A=(A_{1},\ldots,A_{K})$ such that $UU(\beta)=BA$. 
\\
Set $$U(\omega)=\left( U(\omega_{1}),\ldots,U(\omega_{K-1}),-\sum_{k=1}^{K-1}U(\omega_{k})\right),$$  $$U(b)=\left( U(b_{1}),\ldots,U(b_{K})\right)$$ and recall that $$\omega^{\star}=\left(\omega^{\star}_{1},\ldots,\omega^{\star}_{K-1},1-\sum_{k=1}^{K-1}\omega^{\star}_{k}\right).$$
Whenever $R$ is a $K$-dimensional vector, denote $Diag(R)$  the $K\times K$ diagonal matrix having the $R_{k}$'s on the diagonal. 
\\
 Let $P,$ $Q$ and  $\Delta$ be diagonal matrices such that $P = Diag (U(\omega)),$ $Q=Diag(U(b))$ and $\Delta=Diag(\omega^{\star}).$ For $W\in  \mathbb{R}{d},$ set, $D=Diag\left(\langle \beta_1,W \rangle,\ldots,\langle \beta_K,W \rangle\right).$ Then using the fact that $B$ is full rank, (\ref{3}) gives that
 \begin{eqnarray}
&&G_{3} P D+ G_{4} \Delta Q D+ G_{5} +AG_{3} \Delta D + G_{3} \Delta DA^{T}
+\Delta Diag\left(\langle U(\beta_1),\beta_1\rangle,\ldots,\langle U(\beta_K),\beta_K\rangle\right),BA_{K}\rangle)D 
\nonumber\\
&& 
+ G_{3} Diag\left(\langle U(\beta_1),\beta_1\rangle,\ldots,\langle U(\beta_K),\beta_K\rangle\right)=\overline{0}\otimes \overline{0}.\label{3bis}
\end{eqnarray}
Since $U(\beta)=BA,$ then $U(\beta_k)=\sum_{r=1}^K\beta_rA_{rk}=BA_k.$ 
This implies that $$Diag(\langle U(\beta_1),\beta_1\rangle,\ldots,\langle U(\beta_K),\beta_K\rangle)=Diag(\langle BA_1,\beta_1\rangle,\ldots,\langle A_K,\beta_K\rangle)=\tilde{D}.$$
So (\ref{3bis}) can be rewriten as
\begin{eqnarray}
G_{3} P D+ G_{4} \Delta Q D+ G_{5} \Delta \tilde{D} D +AG_{3} \Delta D
 + G_{3} \Delta DA^{T} + G_{3} \Delta \tilde{D}=\overline{0}\otimes \overline{0},\label{3bis1}
\end{eqnarray}
So that for all $W\in  \mathbb{R}{d}$, $W\in  \mathbb{R}{d}$, $AG_{3} \Delta D
 + G_{3} \Delta DA^{T}$ is a diagonal matrix. Since $G_{3} \Delta$ has no zero entries, this proves that, under (H1) and (H3), $A$ is a diagonal matrix. In such a case, $$\tilde{D}=A\tilde{B}\,\,\,\mbox{avec}\,\,\, \tilde{B}=Diag \left(\|\beta_{1}\|^{2},\ldots,\|\beta_{K}\|^{2}\right)$$  
and (\ref{3bis1}) can be rewriten as
\begin{eqnarray}
G_{3} P D+ G_{4} \Delta Q D+ G_{5} \Delta A \tilde{B} D +AG_{3} \Delta D + G_{3} \Delta DA^{T} + G_{3} \Delta A \tilde{B}=\overline{0}\otimes \overline{0}.\label{3ter}
\end{eqnarray}
But by taking, for $k=1,\ldots,K$, $W_{k}$ such that $\beta_{k}^{T}W_{k}=0$, we have $D=0.$ In this case, (\ref{3ter}) is given by $$G_{3} \Delta A \tilde{B}=\overline{0}\otimes \overline{0},$$ 
and using the fact that  we get that $G_{3}$, $\Delta$ and  $\tilde{B}$ have no zero entries we get that $A=0$. This implies that $U(\beta_k)=0,\,\,\,k=1,2,\ldots,K.$ Then using the fact that $B$ is full rank, we conclude from (\ref{1}) and (\ref{2}) that 

\begin{equation}
\label{1bis}
G_{1} P + G_{2} \Delta Q =\overline{0}\otimes \overline{0},
\end{equation}
and
\begin{equation}
\label{2bis}
G_{2} P + G_{3} \Delta Q =\overline{0}\otimes \overline{0}.
\end{equation}
Multipliying (\ref{1bis}) by $G_3$ and (\ref{2bis}) by $G_2,$ we have


\begin{equation}
\label{1bis1}
G_{1}G_3 P + G_{2}G_3 \Delta Q =\overline{0}\otimes \overline{0},
\end{equation}
and
\begin{equation}
\label{2bis1}
G_{2}^2 P + G_2G_{3} \Delta Q =\overline{0}\otimes \overline{0}.
\end{equation}
Taking the difference (\ref{1bis1})-(\ref{2bis1}), we get $$\left( G_{1}G_{3}-G_{2}^{2}\right)P=\overline{0}\otimes \overline{0},$$
and since $G_{1}G_{3}-G_{2}^{2}$ has no zero entries, this leads to $P=0.$ Moreover, since $G_3\Delta$ has no zero entries, this leads also $Q=0$. Thus, under (H1), (H4) and (H5), the matrix $V$ is full rank.

\subsection{Proof of Theorem \ref{theobis1}}
\label{sec:theobis1}
If one knows the law of $(Y,X)$ then the function
$$x \mapsto H(x)=\sum^{K}_{k=1}\omega_k g(\lambda_{k}\langle \mu_k,x \rangle+b_k)
$$ is known on the support of $X$, thus on $ \mathbb{R}^d$.
Then the function
$$
DH(x)= \sum^{K}_{k=1}\omega_k g'(\lambda_{k}\langle \mu_k,x \rangle+b_k)\mu_{k}
$$
is known, 
and if $V\in \mathbb{R}^d$, $\lim_{t\rightarrow +\infty} \| DH(tV)\|=0$ except in case $V$ is  orthogonal to at least one of the $\mu_{k}$'s. 
The set of $V\in \mathbb{R}^d$ such that $\lim_{t\rightarrow +\infty} \| DH(tV)\|\neq 0$ is then $\cup_{k=1}^{K}\langle \mu_{k} \rangle ^{\perp}$, union of disjoint vectorial spaces of dimension  $d-1$, which allows to recover the orthogonal space of $\langle \mu_{k} \rangle ^{\perp}$ for all $k$, thus to recover $K$and all one dimensional spaces $\langle \mu_{k}\rangle$. Since for all $k$, $\omega_{k}g'(b_{k})>0$, this allows to recover the $\mu_{k}$'s. 

\subsection{Proof of Theorem \ref{theobis2}}
\label{sec:theobis2}
Using Theorem \ref{theobis1}, one knows $K$ and the $\mu_{k}$'s. Since the $\mu_{k}$'s are linearly independent, by considering the spaces that  are orthogonal to all  $U_{k}$'s except one, we see that the following functions are known: $h_{1},\ldots,h_{K}$  given for $j=1,\ldots,K$ by:
$$
t \mapsto h_{j}(t)=\omega_j g(\lambda_{j}t+b_j)+ \sum^{K}_{k=1, k\neq j}\omega_k g(b_k).
$$
Then:
\begin{eqnarray*}
h_{j}(0)&=&\sum^{K}_{k=1}\omega_k g(b_k),\\
\lim_{t\rightarrow +\infty} h_{j}(t)&=&\omega_{j} + \sum^{K}_{k=1, k\neq j}\omega_k g(b_k),\\
\lim_{t\rightarrow -\infty} h_{j}(t)&=& \sum^{K}_{k=1, k\neq j}\omega_k g(b_k).\\
\end{eqnarray*}
This allows to recover $\omega_{j}$ and $g(b_{j})$ for $j=1,\ldots,K$. Thus the functions
$$
t \mapsto \ell_{j}(t)= g(\lambda_{j}t+b_j)
$$
are known.
Since $g(0)=\ell_{j}(-b_{j}/\lambda_{j})$ and $g(1)=\ell_{j}((1-b_{j})/\lambda_{j})$ are fixed, one can find $\lambda_{j}$ and $b_{j}$, and then the function $g$.

\subsection{Proof of Theorem \ref{theobis3}}
\label{sec:theobis3}
Using Theorem \ref{theobis1} applied to the distributions of $Y$ conditional to  $X$ and $Z=z$ for all $z\in \{z_{1},\ldots,z_{m}\}$, the knowledge of $ P_{g,\omega,\beta,\gamma,b}$  allows to recover $K$, $g$, $\omega$, $\beta$, and $A_{k}=(a_{k,i})_{1\leq i \leq m}$, $k=1,\ldots,K$, with
$$
a_{k,i}=b_{k}+ \langle \gamma_{k},z_{i}\rangle.
$$
We then know for all $k$
$$
A_{k}=\left(\begin{array}{ll}1 & z_{1}^{T}\\
1 & z_{2}^{T}\\
\vdots & \vdots \\
1 & z_{m}^{T}
\end{array}\right)
\left(\begin{array}{l}b_{k}\\
\gamma_{k}
\end{array}\right)
$$
which allows to recover the $b_{k}$'s and $\gamma_{k}$'s when (S4) holds.

\subsection{Proof of Theorem \ref{theobis4}}
\label{sec:theobis4}
If one knows the law of $Y$, then, for all fixed $z\in  \{z_{1},\ldots,z_{m}\}$, one knows the function $H:(\mathbb{R}^{d})^{p}\rightarrow (0,1)^{p}$ given by
$$
H(x_{1},\ldots,x_{p})=\sum^{K}_{k=1}\omega_k 
\left(g(\langle \beta_k,x_{j} \rangle+  \tilde{b}_k (z))\right)_{1\leq j \leq p}
$$
with $\tilde{b}_k (z)=b_{k}+ \langle \gamma_{k},z_{i}\rangle.$
Let us firts prove that for all $z$, the functions $g(\langle \beta_k,\cdot \rangle+\tilde{b}_k (z))$ are linearly independent. Indeed, if $\alpha_{1},\ldots,\alpha_{K}$are such that for all $x\in \mathbb{R}^d$,
$$
\sum_{k=1}^{K} \alpha_{k}g(\langle \beta_k,x \rangle+\tilde{b}_k (z))=0,
$$
then by taking the derivative, for all  $x\in \mathbb{R}^d$,
$$
\sum_{k=1}^{K} \alpha_{k}g'(\langle \beta_k,x \rangle+\tilde{b}_k (z))\beta_{k}=0.
$$
Since (S2) holds, 
there exists $V \in \langle \beta_{k} \rangle^{\perp}$ such that $V \notin \langle \beta_{j} \rangle^{\perp}$, $j\neq k$. Then taking $x=tV$ and $t$ tending to infiny, we get that $\alpha_{k}g'(\tilde{b}_k (z))\beta_{k}=0$, and then $\alpha_{k}=0$. \\
Now, following the spectral method of proof  developed in \cite{AK2014} to prove that multidimensional mixtures are identifiable (see also \cite{MR3889980} Chapter 14) , we see that the knowledge of $H$ allows to recover $K$, the $\omega_{k}$'s and, for all $z$, the functions $g(\langle \beta_k,\cdot \rangle+\tilde{b}_k (z))$.\\
Then, if one knows the function $x\mapsto g(\lambda_{k}\langle \mu_{k},x \rangle + \tilde{b}_k (z))$ one can recover $\mu_{k}$ by taking the derivative, then $g$ and the $\tilde{b}_k (z))$'s as in the proof of Theorem \ref{theobis2} then the $\gamma_{k}$'s and the $b_{k}$'s  as in the proof of Theorem \ref{theobis3}.
